\magnification=1100
\baselineskip=12pt
\hsize12cm
\vsize18cm
\font\twelverm=cmr12
\font\twelvei=cmmi12
\font\twelvesy=cmsy10
\font\twelvebf=cmbx12
\font\twelvett=cmtt12
\font\twelveit=cmti12
\font\twelvesl=cmsl12

\font\ninerm=cmr9
\font\ninei=cmmi9
\font\ninesy=cmsy9
\font\ninebf=cmbx9
\font\ninett=cmtt9
\font\nineit=cmti9
\font\ninesl=cmsl9

\font\eightrm=cmr8
\font\eighti=cmmi8
\font\eightsy=cmsy8
\font\eightbf=cmbx8
\font\eighttt=cmtt8
\font\eightit=cmti8
\font\eightsl=cmsl8

\font\sixrm=cmr6
\font\sixi=cmmi6
\font\sixsy=cmsy6
\font\sixbf=cmbx6

\catcode`@=11 
\newskip\ttglue
\def\twelvepoint{\def\rm{\fam0\twelverm}
\textfont0=\twelverm  \scriptfont0=\ninerm  
\scriptscriptfont0=\sevenrm
\textfont1=\twelvei  \scriptfont1=\ninei  \scriptscriptfont1=\seveni
\textfont2=\twelvesy  \scriptfont2=\ninesy  
\scriptscriptfont2=\sevensy
\textfont3=\tenex  \scriptfont3=\tenex  \scriptscriptfont3=\tenex
\textfont\itfam=\twelveit  \def\it{\fam\itfam\twelveit}%
\textfont\slfam=\twelvesl  \def\sl{\fam\slfam\twelvesl}%
\textfont\ttfam=\twelvett  \def\tt{\fam\ttfam\twelvett}%
\textfont\bffam=\twelvebf  \scriptfont\bffam=\ninebf
\scriptscriptfont\bffam=\sevenbf  \def\bf{\fam\bffam\twelvebf}%
\tt  \ttglue=.5em plus.25em minus.15em
\normalbaselineskip=15pt
\setbox\strutbox=\hbox{\vrule height10pt depth5pt width0pt}%
\let\sc=\tenrm  \let\big=\twelvebig  \normalbaselines\rm}

\def\tenpoint{\def\rm{\fam0\tenrm}
\textfont0=\tenrm  \scriptfont0=\sevenrm  \scriptscriptfont0=\fiverm
\textfont1=\teni  \scriptfont1=\seveni  \scriptscriptfont1=\fivei
\textfont2=\tensy  \scriptfont2=\sevensy  \scriptscriptfont2=\fivesy
\textfont3=\tenex  \scriptfont3=\tenex  \scriptscriptfont3=\tenex
\textfont\itfam=\tenit  \def\it{\fam\itfam\tenit}%
\textfont\slfam=\tensl  \def\sl{\fam\slfam\tensl}%
\textfont\ttfam=\tentt  \def\tt{\fam\ttfam\tentt}%
\textfont\bffam=\tenbf  \scriptfont\bffam=\sevenbf
\scriptscriptfont\bffam=\fivebf  \def\bf{\fam\bffam\tenbf}%
\tt  \ttglue=.5em plus.25em minus.15em
\normalbaselineskip=12pt
\setbox\strutbox=\hbox{\vrule height8.5pt depth3.5pt width0pt}%
\let\sc=\eightrm  \let\big=\tenbig  \normalbaselines\rm}

\def\ninepoint{\def\rm{\fam0\ninerm}
\textfont0=\ninerm  \scriptfont0=\sixrm  \scriptscriptfont0=\fiverm
\textfont1=\ninei  \scriptfont1=\sixi  \scriptscriptfont1=\fivei
\textfont2=\ninesy  \scriptfont2=\sixsy  \scriptscriptfont2=\fivesy
\textfont3=\tenex  \scriptfont3=\tenex  \scriptscriptfont3=\tenex
\textfont\itfam=\nineit  \def\it{\fam\itfam\nineit}%
\textfont\slfam=\ninesl  \def\sl{\fam\slfam\ninesl}%
\textfont\ttfam=\ninett  \def\tt{\fam\ttfam\ninett}%
\textfont\bffam=\ninebf  \scriptfont\bffam=\sixbf
\scriptscriptfont\bffam=\fivebf  \def\bf{\fam\bffam\ninebf}%
\tt  \ttglue=.5em plus.25em minus.15em
\normalbaselineskip=11pt
\setbox\strutbox=\hbox{\vrule height8pt depth3pt width0pt}%
\let\sc=\sevenrm  \let\big=\ninebig  \normalbaselines\rm}

\def\eightpoint{\def\rm{\fam0\eightrm}
\textfont0=\eightrm  \scriptfont0=\sixrm  \scriptscriptfont0=\fiverm
\textfont1=\eighti  \scriptfont1=\sixi  \scriptscriptfont1=\fivei
\textfont2=\eightsy  \scriptfont2=\sixsy  \scriptscriptfont2=\fivesy
\textfont3=\tenex  \scriptfont3=\tenex  \scriptscriptfont3=\tenex
\textfont\itfam=\eightit  \def\it{\fam\itfam\eightit}%
\textfont\slfam=\eightsl  \def\sl{\fam\slfam\eightsl}%
\textfont\ttfam=\eighttt  \def\tt{\fam\ttfam\eighttt}%
\textfont\bffam=\eightbf  \scriptfont\bffam=\sixbf
\scriptscriptfont\bffam=\fivebf  \def\bf{\fam\bffam\eightbf}%
\tt  \ttglue=.5em plus.25em minus.15em
\normalbaselineskip=9pt
\setbox\strutbox=\hbox{\vrule height7pt depth2pt width0pt}%
\let\sc=\sixrm  \let\big=\eightbig  \normalbaselines\rm}

\def\twelvebig#1{{\hbox{$\textfont0=\twelverm\textfont2=\twelvesy
	\left#1\vbox to10pt{}\right.\n@space$}}}
\def\tenbig#1{{\hbox{$\left#1\vbox to8.5pt{}\right.\n@space$}}}
\def\ninebig#1{{\hbox{$\textfont0=\tenrm\textfont2=\tensy
	\left#1\vbox to7.25pt{}\right.\n@space$}}}
\def\eightbig#1{{\hbox{$\textfont0=\ninerm\textfont2=\ninesy
	\left#1\vbox to6.5pt{}\right.\n@space$}}}
 
\def\tit{\bigskip}

\def\Pt{{\bf P}^3}
\def\v{{\cal V}}

\def\o{{\cal O}}
\def\i{{\cal I}}

\def\e{{\cal E}}
\def\l{{\cal L}}
\topinsert 
\endinsert
\font\big=cmbx10 scaled \magstep2

\centerline{\big \hbox {Gonality, Clifford index and multisecants}}
\vskip 0.8truecm
\centerline {Ph. Ellia $^*$ - D. Franco $^{**}$ \footnote {$^1$}
{Partially supported by MURST and Ferrara Univ. in the framework
of the project:
"Geometria algbrica, algebra commutativa e aspetti computazionali"}}
\medskip
{\eightpoint 
\centerline {Dipartimento di Matematica, Universit\`a di Ferrara}
\centerline {via Machiavelli 35 - 44100 Ferrara, Italy}
\centerline {$^{*}$ e-mail: phe@dns.unife.it}
\centerline {$^{**}$ e-mail: frv@dns.unife.it}}
\bigskip

\centerline {\bf Introduction }
\tit
The gonality of a projective curve $C$ is the minimal degree
of a morphism $f: C\to {\bf P}^1$. It is a classical invariant 
which has been refined by the introduction of the Clifford index:
$$
Cliff(C):=min\{Cliff(\l ): \hskip1mm \l \in Pic(C) , 
\hskip1mm h^0(\l )\geq 2,
\hskip1mm h^1(\l )\geq 2 \}
$$
where $Cliff(\l ):=deg(\l ) -2(h^0(\l )-1)$
(see for example [CM] and [ELMS]).
\par  
In the first section of this paper we consider
subcanonical curves in $\Pt $.
Let  $C\subset \Pt$ 
and let  $\Gamma $ be a point set computing the gonality of $C$.
If $l\geq 2$ represents the maximum degree of a 
zero-dimensional subscheme
of $C$ which is contained in a line, 
then $d:=Gon(C)=d(\Gamma)\leq d(C)-l$.
If $Gon(C)= d(C)-l$ we will say that {\it Gon(C) is computable
by multisecants }. In [B] Basili proved that if $C$ is a complete intersection
then $Gon(C)$ is computable by multisecants. 
Furthermore in the same paper, Basili computes the Clifford index
of complete intersections.
In this paper we generalize these results
to most subcanonical curves in $\Pt $ (see Theorem 1.11):
\tit
{\bf Theorem. }{\it Let $E$ be a rank two vector bundle in $\Pt$.
If $t>>0$ and if $C$ is a curve which is the zero locus of a section of $E(t)$
then $Gon(C)$ is computable by multisecants.
Moreover, 
if $C$ is  not bielliptic,
then either $Cliff(C)=Gon(C)-3$ or $Cliff(C)=Gon(C)-2$. 
Furthermore, the following conditions are equivalent:
\item{1)}{$Cliff(C)=Gon(C)-3$;}
\item{2)}{$Cliff(C)=Cliff(\o _C(1))=d(C)-6$;} 
\item{3)}{$C$ does not have four-secant lines (i.e. $l=3$).}}
\tit
Our approach is completely different from Basili's one and relies on 
vector bundle techniques: Bogomolov's unstability theorem [La] and
Tyurin's work [T].

In section 2 we consider a natural problem arising already from Basili's work:
the stratification by multisecants
of the Hilbert scheme of complete intersections.
This problem is interesting by itself, not only for complete intersections,
and we consider it in two extremal and opposite situations:
complete intersections and rational curves.
In both cases we prove that the locus of  curves with a k-secant
line is irreducible and of the expected dimension except when this
cannot be true for trivial reasons (see Remark 2.7).

To conclude let us suggest two directions to extend the results of this paper:

investigate the stratification by multisecants for other Hilbert schemes;

determine further classes of curves with gonality computable by multisecants.
\par\noindent
Actually this paper originated by a suggestion of Peskine that Basili's
result should extend to projectively normal curves (notice that
a complete intersection is both subcanonical and projectively normal).
\tit
{\bf Acknowledgements. } We thank Ch. Peskine for drawing our attention
on Basili's paper.
\tit
\centerline{\bf Gonality and Clifford index }
\tit
Consider a smooth subcanonical curve $C\subset \Pt$ ($\omega _C\simeq \o _C(\alpha)$)
and a point set $\Gamma $ computing the gonality of $C$.
If $l\geq 2$ represents the maximum degree of a zero-dimensional subscheme
of $C$ which is contained in a line, then $d:=Gon(C)=d(\Gamma)\leq d(C)-l$.
\tit
{\bf Definition 1.1. }{\it If Gon(C)=d(C)-l we will say that Gon(C) is computable by
multisecants }
\tit
In this section we are going to prove
the equality in the last formula for several classes of subcanonical curves. 
\par
We start with a technical lemma.
\tit
{\bf Lemma 1.2.}{ \it 
\par
Suppose $h^1(\i _C(\alpha))=0$.
Then $h^0(\o _C(\Gamma))=h^1(\i _{\Gamma,\Pt}(\alpha))+1$. }
\tit
{\bf Proof.} Consider the sequence
$$
0\to \i _C(\alpha)\to \i _{\Gamma}(\alpha)\to \i _{\Gamma,C}(\alpha)\to 0
$$
and look at the cohomology sequence.
Since $h^1(\i _C(\alpha))=0$ and 
$H^2(\i _C(\alpha))\simeq H^1(\o _C(\alpha))\simeq H^0(\o _C)^*\simeq {\bf C}$
we find 
$$
0\to H^1(\i _{\Gamma}(\alpha))\to H^1(\i _{\Gamma,C}(\alpha))\to {\bf C}\to 0
$$
and the thesis follows taking account of 
$H^1(\i _{\Gamma,C}(\alpha))\simeq H^0(\o _C(\Gamma))$.
\tit
{\bf Remarks 1.3. } a) As in the lemma above, suppose $h^1(\i _C(\alpha ))=0$.
Since $h^1(\i _{\Gamma,\Pt}(\alpha))=h^1(\i _{\Gamma,F }(\alpha))$
for any surface $F\subset \Pt$, in order that  $\Gamma $ contributes to the gonality
of $C$ it is necessary and sufficient that $\Gamma $ fails to impose independent
conditions to the linear series $\alpha H$ on $F$.
From the very definition of $\Gamma $, that is from the condition 
$Gon(C)=d(\Gamma )$, it is clear that $h^1(\i _{\Gamma ',F}(\alpha))=0$
for each proper subset $\Gamma '\subset \Gamma $.
This implies that the set $\Gamma $ is $\alpha H$ stable for any 
smooth surface containing it (see [T] and [P]) hence $\Gamma $ is the zero locus
of a section of a rank two vector bundle $\e $ on $F$ sitting in a sequence 
like the following [T]:
$$
0\to \o _{F} \to \e \to \i _{\Gamma,F}(\alpha-f+4)\to 0.
$$
Furthermore, the fact that $\Gamma $ computes the gonality of $C$
implies that 
\par\noindent
$h^0(\o _C(\Gamma ))=2$ hence, by Lemma 1.1, 
$h^1(\i _{\Gamma ,F}(\alpha))=1$ and the vector bundle $\e $ 
is uniquely determined.
\par
b) Keep the notations of the previous remark and
consider a $pH$ stable divisor $\Gamma \subset C \subset F$ 
such that $h^1(\i _{\Gamma ,F}(p))\geq 1$. 
Then Tyurin's technique produces a vetor bundle $\v $ 
of rank $h^1(\i _{\Gamma ,F}(p))+1$ whose general rank-two quotient $\e $ is a vector
bundle embedded in a sequence like in the previous remark
(of course now the bundle $\e $ is not uniquely determined):
$$
0\to \o _{F} \to \e \to \i _{\Gamma ,F}(p-f+4)\to 0.
$$
Suppose now $h^1(\i _C(\alpha ))=0$.
The construction just outlined applies in particular to any divisor $\Gamma $ 
computing the Clifford index of $C$ with $p=\alpha $. 
Indeed it is easy to see that
$\Gamma $ is $\alpha H$-stable: 
denotes by $\Gamma '\subset \Gamma $  the stable part
of $\Gamma $ (see [T] for the definition of stable part of a zero dimensional
scheme in a smooth surface). We have 
$h^0(\o _C(\Gamma '))=h^0(\o _C(\Gamma ))$ hence
$Cliff(\Gamma ')\leq Cliff (\Gamma )$ with equality iff $\Gamma '=\Gamma $.
\par
c) Set $s(C):=min\{n: \hskip2mm h^0(\i _C(n)\not =0 \}$ and consider a smooth surface
$F$ containing $C$ whose degree $f$ is bigger than $s(C)$. Suppose that there exists 
a smooth surface $T$, with $deg(T)<f$, such that $T\cap F=C\cup D$ with $D$ integral. 
A  theorem of Lopez 
(see [L]) says that the Picard group of the general surface of degree $f$ containing $C$, has rank two and is generated by 
the hyperplane and $D$.
\tit
{\bf  Notation. } {\it In the sequel we will make the following
assumption:
\item{($\circ $)}{ there exist smooth surfaces $T$,$F'$,$f=deg(F')>deg(T)$ 
containing $C$ such that $T\cap F'=C\cup D$ where $D$ is smooth.}
\par
\noindent
By Remark $c)$ above  for the general surface $F$ of degree $f$ containing $C$
we have $Pic(F)=<H,D>$.}
\tit
{\bf Theorem 1.4.}{ \it Let $C$ be a subcanonical ($\omega _C\simeq \o _C(\alpha)$),
curve. Assume condition $(\circ )$ is verified. 
Suppose that $\Gamma\subset C$ 
is a $pH$-stable ($p\leq \alpha $) divisor (in $F$) such that 
$d=d(\Gamma )<d(C)$ and $h^1(\i _{\Gamma ,F }(p))\geq 1$.
Assume the following conditions are verified
\item{a)}{ $h^1(\i _C(1))=0$; }
\item{b)}{ $f<\alpha +4$;}
\item{c)}{ there exists a positive integer $s$ such that $p-f+4>s+{d\over sf}$;}
\item{d)}{ $d(C)\leq 2(p-f+2)f $.} 
\par
Then $\Gamma $ is a planar set and $d>(p-f+3)f$. }
\tit
{\bf Remark 1.5.} Note that condition b) implies that $C^2> 0$ in $F$ hence
$C\cdot E \geq 0$ for any effective divisor $E$.
\tit
{\bf Proof.} By Lemma 1.2 and Remarks 1.3 we find
 a vector bundle $\e $ on $F$ admitting a section vanishing on $\Gamma $.
The discriminant of $\e $ is $\Delta(\e )=(p-f+4)^2f-4d$.
The hypothesis $c)$ above implies  $\Delta(\e )>(s\sqrt f -{d\over s\sqrt f})^2$
hence, by the Bogomolov's unstability theorem (see [La] \S 4),
there exists a divisor $Y\subset F$ ($0\to \o _F(Y)\to \e$)
such that
\item{1)}{ $2Y-(p-f+4)H\in N^+(F)$ (the positive cone of $F$); }
\item{2)}{ $[2Y-(p-f+4)H]^2>\Delta(\e )$. }
\par\noindent
Notice that the last condition is equivalent
to $d>Y\cdot[(p-f+4)H-Y]$.
If we set $X:=(p-f+4)H-Y$ we get
\item{i)}{ $Y-X\in N^+(F)$ ; }
\item{ii)}{ $d>X \cdot Y$. }
\par\noindent 
We claim that condition $1)$ implies that 
$H^0(\o _F(-Y))=0$. Indeed, $(2Y-kH)\in N^+$ implies
$(2Y-kH)\cdot H>0$ hence $Y\cdot H>0$ and the claim follows. 
Hence the composite morphism 
$$
\o _F(Y) \to \e \to \i _{\Gamma,F}(p-f+4)
$$
does not vanish and $h^0(\i _{\Gamma,F}(p-f+4-Y))\not = 0$.
So we may suppose $X$ to be an effective curve containing $\Gamma $.
Moreover, combining  Remark 1.5 with the fact that some multiple of 
$Y-X$ is effective (see [H] Corollary V 1.8), we find
$$ 
C\cdot Y \geq C\cdot X \eqno (*)
$$
\par
By hypothesis, $Pic(F)=<H,D>$ where $F\cap T=C\cup D$; furthermore  the exact sequence of liaison: 
$$
0\to \i _{F\cap T}(t+f-\alpha -4) \to \i _D(t+f-\alpha -4)
\to \omega _C(-\alpha ) \to 0
$$
shows that there exists an effective divisor $R$ on $F$ 
such that $D+R=(t+f-\alpha -4)H$.
From the exact sequence above $R\cap C=\emptyset $
(because $\omega _C(-\alpha )\simeq \o _C$). Furthermore
$Pic(F)=<H,R>$ and $\omega _R\simeq \o _R(g-4)$ with $g=2f-\alpha -4$.
We have
$Y\sim \beta H+\delta R$ hence 
$X\sim (p-f+4-\beta )H-\delta R$ and
ii) says that 
$$
d>X \cdot Y =\beta (p-f+4-\beta )f+\delta d(R)(p-f+4-2\beta )
- \delta ^2R \cdot R  \eqno(**).
$$
The adjunction formula provides to the selfintersection of $R$
$$
R \cdot R=(g-f)d(R)
$$
hence 
$$
X \cdot Y =\beta (p-f+4-\beta )f+\delta d(R)((p-f+4-2\beta )
- \delta (g-f)) .
$$
By  Lemma  1.7 stated below we find $X\cdot Y\geq \beta (p-f+4-\beta )f$
and Lemma 1.6 combined with $(**)$ implies
$d(C)>d>\beta (p-f+4-\beta )f$ with $0<\beta<p-f+4$.
The last inequality,  the hypothesis
$d(C)\leq 2(p-f+2)f$ 
and the fact that  the maximum of $\phi (\beta ):= \beta (p-f+4-\beta )$
is reached for $\beta ={p-f+4 \over \beta }$
imply that either
$\beta =1$ or $p-f+4-\beta =1$.
We claim that the second item is the right one.
Indeed, $Y\sim \beta H+\delta R$, $X\sim (p-f+4 -\beta )H-\delta R$ and 
$(*)$ implies $ C\cdot Y =\beta d(C)\geq C\cdot X=(p-f+4 -\beta )d(C)$.
Hence $X$ is linearly equivalent to $H-\delta R$.
To prove that $\Gamma $ is planar it suffices to combine $\Gamma \subset X\cap C$ and 
$X\cdot C\sim _C H\cdot C$ with the fact that $C$ is linearly normal.
Finally, the thesis $d>(p-f+3)f$ follows directly from $(**)$ with $\beta =1$.
\tit
{\bf Lemma 1.6. }{\it $0<\beta <p-f+4$. }
\tit
{\bf Proof. } 
\par $\beta <p-f+4$: 
Since $\Gamma\subset X\cap C$ we find 
$C\cdot X=C\cdot [(p-f+4-\beta )H-\delta R]\geq d$.
But $C\cap R= \emptyset $ hence $p-f+4>\beta $.  
\par $0<\beta $:  $\beta d(C)=(\beta H+\delta R)\cdot C=Y\cdot C\geq X\cdot C\geq d>0$ 
(as above, recall that $R\cdot C=0$ and that $\Gamma \subset X\cap C$). 
\tit 
 {\bf Lemma 1.7. } {\it We have
$\delta ((p-f+4-2\beta )+ \delta (f-g))\geq 0$. }
\tit
{\bf Proof. }
We proceed in two steps depending on the 
sign of $\delta $.
\par $\delta\geq 0$: by the last lemma $\beta <p-f+4 $ hence
$p -f+4-2\beta >-(p-f+4)\geq -(f-g)$ (recall that $p\leq \alpha $). The claim follows since
$\delta ((p-f+4-2\beta )+ \delta (f-g))
\geq \delta ( \delta (f-g)-(f-g))\geq 0 $ (note that hypothesis $b)$ implies 
$f-g>0$). 
\par $\delta <0$:  by $(*)$ we have
 $(Y-X)\cdot C\geq 0$ hence $2\beta -p-4+f\geq 0$.
The lemma follows. 
\tit
{\bf Corollary 1.8. }{ \it Let $C$ be a subcanonical ($\omega _C\simeq \o _C(\alpha)$),
curve
and $\Gamma\subset C$ a  divisor such that $d:= d(\Gamma )=Gon(C)$.
Suppose condition $(\circ )$ is verified.
Additionally, assume the  conditions $a)$- $d)$ of Theorem 1.4 with $\alpha =p$
and suppose $h^1(\i _C(\alpha ))=0$.
\par
Then $\Gamma $ is a planar set hence $Gon(C)=d(C)-l$ i.e. $Gon(C)$ is
computable by multisecants. }
\tit
{\bf Proof. } It follows combining Lemma 1.2, Remark 1.3 $a)$ and Theorem 1.4.
\tit
{\bf Corollary  1.9. }{ \it Let $C$ be a subcanonical ($\omega _C\simeq \o _C(\alpha)$),
curve
and $\Gamma\subset C$ a  divisor.
Supose condition $(\circ )$ holds and assume the 
following conditions are verified (as above, $p\leq \alpha$)
\item{a)}{ $h^1(\i _C(1))=0$; }
\item{b)}{ $f<\alpha +4$;}
\item{c)}{ there exists a positive integer $s$ such that $p-f+4>s+{d(\Gamma ) \over sf}$;}
\item{d)}{ $d(C)\leq 2(p-f+2)f $.} 
\par
If $\Gamma _1\subset \Gamma $ is a subdivisor such that 
$d(\Gamma _1)\leq (p-f+3)f$ then 
\par\noindent $h^1(\i  _{\Gamma _1,F }(p))=0$. }
\tit
{\bf Proof. } If $h^1(\i  _{\Gamma _1,F }(p))\not =0$, applying
Theorem 1.4 to the stable part of $\Gamma _1 $ we should get $d> (p-f+3)f$ 
(see [T] for the definition of stable part).
\tit
{\bf Theorem 1.10. }{ \it Let $C$ be a subcanonical 
($\omega _C\simeq \o _C(\alpha)$, $\alpha\geq 4$)
curve
and $\Gamma\subset C$ a  divisor computing the Clifford index of $C$.
Assume  $(\circ )$.
Suppose further that $C$ is neither hyperelliptic nor bielliptic.
Set $d:=d(\Gamma )$.
Assume the following conditions are verified
\item{a)}{ $h^1(\i _C(\alpha ))=h^1(\i _C(1))=0$; }
\item{b)}{ $f<\alpha +4$;}
\item{c)}{ there exists a positive integer $s$ such that $\alpha -f+3>s+{d\over sf}$;}
\item{d)}{ $d(C)\leq 2(\alpha -f+1)f $.} 
\par
Then either $Cliff(C)=Gon(C)-3$ or $Cliff(C)=Gon(C)-2$. 
Moreover, the following conditions are equivalent:
\item{1)}{$Cliff(C)=Gon(C)-3$;}
\item{2)}{$Cliff(C)=Cliff(\o _C(1))=d(C)-6$;} 
\item{3)}{$C$ does not have four-secant lines (i.e. $l=3$).}  }
\tit
{\bf Proof. } 
Note that the numerical conditions are just the same of Theorem 1.4 for 
$p=\alpha -1$ and that they are {\it a fortiori } verified for 
$p=\alpha $. Hence we can apply Corollary 1.8 and $Gon(C)=d(C)-l$.
\par
Arguing as in Theoreme 4.3 of [B] it easy to show that
either $Cliff(C)=Gon(C)-2$ or $Cliff(C)=Gon(C)-3$.
Let us suppose that $Cliff(C)=Gon(C)-3$. Then (see again Theoreme 4.3 of [B])
there exists a subdivisor $\Gamma '\subset \Gamma $ such that 
$Gon(C)=d(\Gamma ')=d(C)-l$. We know that there exists a plane $H$ containing
$\Gamma '$. Since $\alpha \geq 4$ we find 
$g(C)={\alpha d(C)+2 \over 2}> 2(d(C)-l-3)+5=2Cliff(C)+5$
hence, by [CM] Cor. 3.2.5, 
$$
d(\Gamma )\leq {3\over 2}(Cliff(C)+2). \eqno (+)
$$
\par
Let $\Gamma _1 $ be the residual set to $\Gamma \cap H$ inside $\Gamma $.
By $(+)$ the degree of $\Gamma _1$ is bounded by 
${1\over 2}d(C) -{l+1\over 2}$. The hypothesis 
$d(C)\leq 2(\alpha -f+1)f$ implies $d(\Gamma _1)<(\alpha -f+2)f$ hence,
by Corollary 1.9 $h^1(\i _{\Gamma _1,F}(\alpha -1))=0$.
\par
By the following sequence
$$
0\to \i _{\Gamma _1,F}(\alpha -1)\to  \i _{\Gamma ,F}(\alpha ) \to 
\i _{\Gamma \cap H,H\cap F}(\alpha )\to 0
$$
we have $h^0(\o _C(\Gamma ))=h^0(\o _C(\Gamma \cap H))$ (recall Lemma 1.2)
hence $\Gamma \subset H$ because $\Gamma $ computes $Cliff(C)$.
Note that the assumption $Cliff(C)=Gon(C)-3$ implies that 
$h^0(\o _C(\Gamma ))\geq 3$ (otherwise $\Gamma $ would compute the gonality of $C$
and $Cliff(C)=Gon(C)-4$).
On the other hand, $\Gamma \subset H$ implies $h^0(\o _C(\Gamma ))\leq 4$ with
equality if and only if $\o _C(\Gamma )\simeq \o _C(1)$.
Combining the last inequalities with $l\geq 3$ we see that the hypothesis
$Cliff(\Gamma )=d(\Gamma )-2(h^0(\o _C(\Gamma ))-1)=Gon(C)-3=d(C)-l-3$
implies $\o _C(\Gamma )\simeq \o _C(1)$ and $l=3$.
\par 
To conclude the proof it suffices to  show that $l=3$ implies
$Cliff(C)=Gon(C)-3$ and $\Gamma=\o _C(1)$ which follows at once by 
$Cliff(C)\geq Gon(C)-3=d(C)-6 =Cliff(\o _C(1))$.
\tit
{\bf Theorem 1.11. }{\it Let $C\subset \Pt$ be  
the zero locus of a section of an high twist of a rank
two vector bundle of $\Pt$.
Then we have $Gon(C)=d(C)-l$.
Suppose $C$ is  not bielliptic.
Then either $Cliff(C)=Gon(C)-3$ or $Cliff(C)=Gon(C)-2$. 
Moreover, the following conditions are equivalent:
\item{1)}{$Cliff(C)=Gon(C)-3$;}
\item{2)}{$Cliff(C)=Cliff(\o _C(1))=d(C)-6$;} 
\item{3)}{$C$ does not have four-secant lines (i.e. $l=3$).}}
\tit
{\bf Proof.}
Let $C$ be a curve coming
as the zero locus of an high order twist of a normalized vector bundle 
$E$ over $\Pt $:
$$
0\to \o \to E(t) \to \i _C(2t+c_1) \to 0 \hskip3mm t>>0
$$
(where either $c_1=0$ or $c_1=-1$ ).
\par
The  assumption a) of Theorem 1.4 is satisfied:
\par
$h^1(\i _C(\alpha))=h^1(\i _C(2t+c_1-4))=h^1(E(t-4))=0$ if $t>>o$.
\par
$h^1(\i _C(1))=h^1(E (1-t-c_1))=0$ if $t>>0$.
\par\noindent
Let $r$ be such that $E(r)$ is globally generated, then also $\i _C(r+c_1+t)$ 
is globally generated and condition $(\circ)$ is
verified for the general surface in $H^0(\i _C(r+c_1+t+1))$. Furthermore, we have 
$r+c_1+t+1>s(C)$.
Now we apply Theorem 1.10 with $f=r+c_1+t+1$ and we notice that 
if $\omega _C\sim \o _C(\alpha)$ with $\alpha >0$ then $C$ is not 
hyperelliptic.
We verify the numerical conditions for $f=r+c_1+t+1$:
\item{b):}{ $r+t+1<2t$; satisfied for $t>r$.}
\item{c):}{ $t-r-2>s+{d\over s(r+c_1+t+1)}$; since $d<d(C)=c_2+c_1t+t^2$, it is enough to 
check $t-r-2\geq s+{c_2+c_1t+t^2\over s(r+c_1+t+1)}$ which is verified by $s=2$ and $t>>0$.}
\item{d):}{ $c_2+c_1t+t^2\leq 2(t-r-4)(r+c_1+t+1)=2[t^2+(c_1-3)t-(r+c_1+1)(r+4)]$; it is 
verified for $t>>0$.}
\tit
{\bf Examples. }

i) In the case of a decomposed vector bundle we can make
precise computations and recover Basili's results except
for complete intersections of type: $(3,a)$, $(4,4)$ and $(5,5)$.

ii) In the same vein if  $C$ be the zero locus of a section
of $N(t)$ ($t\geq 7$) where $N$ is a normalized Null-Correlation 
bundle of $\Pt $
we see that the argument of the previous proof  applies hence 
the conclusion of Theorem 1.11 holds.
\tit
\centerline{\bf Multisecants to space curves.}
\tit
The results of the previous section naturally introduce the following
question: for which $k$ does there exists a smooth complete intersection
curve of type $(a,b)$ with a ''maximal'' $k$-secant line (i.e. the curve has
no $l-$secant line with $l>k$)? More generally one could ask for a
description of the locus of complete intersections with a $k-$secant line.
This natural problem is of interest by itself, and not only for complete
intersections curves. We will consider two extremal, and opposite
situations: complete intersections and rational curves.
\tit
{\bf Generalities}
\tit
{\bf Notations 2.1.} We denote by $H(d,g)$ the closure in $Hilb(\Pt )$ of
the set of smooth, connected curves of degree $d$, genus $g$. By $H^s(d,g)$
we will denote the open subset of $H(d,g)$ parametrizing smooth curves.
Similarly if $H$ is an irreducible component of $H(d,g)$, $H^s$ will denote
the open subset of $H$ corresponding to smooth curves.

We will denote by $Al^k$ the closed subscheme of $Hilb^k\Pt$
parametrizing zero-dimensional subschemes of length $k$ which are contained
in a line (i.e. which are ''aligned''). We recall that, for $k>1$, $Al^k$ is
smooth, irreducible, of dimension $4+k$ (consider the natural map $%
Al^k\rightarrow Gr(1,3)$).

If $H$ is an irreducible component of $H(d,g)$ we have the incidence variety 
$\i _k(H)$ (or, if no confusion can arise $\i _k$): $%
\i _k:=\left\{ (Z,C)\in Al^k\times H\mid Z\subset C\right\} ;$
associated to $\i _k$ we have the diagram:
$$
\matrix{
\i _k & \buildrel p_k \over \rightarrow  & H \cr 
 \downarrow q_k &  &  \cr 
Al^k &  & \cr }
$$
where $p_k,q_k$ denote the natural projections.

In this situation we define $\i _k^s$ as $p_k^{-1}(H^s)$. Clearly $%
\i _k^s$ is open in $\i _k.$
\tit
{\bf Definition 2.2.} {\it The locus of curves of $H$ with a $k-$secant line is 
$H_k:=p_k(\i _k) $. 
The locus of smooth curves of $H$ with a $k-$secant line is 
$H_k^s:=H_k\cap H^s=p_k(\i _k^s)$. }
\tit
{\bf Remarks 2.3.}
\item{(i)}{ A line $L$ is a $k$-secant to the smooth connected curve $C$ (of degree $%
>1)$ if $length(C\cap L)\geq k;L$ is a 
{\it proper $k-$secant} if 
$length(C\cap L)=k,L$ is a 
{\it maximal secant} to $C$ if for $l>k,C$
has no $l-$secant lines.}
\item{(ii)}{ It may happen that $\i _k$ is empty. It may also happen that $
\i _k\neq \emptyset $ while $\i _k^s$ is empty. }
\item{(iii)}{ To pass through one point impose two conditions to curves in $H$, so
we may expect, in general, the general fiber of $q_k$ to be of dimension $h-2k$ 
($h=\dim (H)$), hence we may expect $\dim (\i _k)=h-(k-4).$
Also, in general, and if $k>3,$ we may expect the general fiber of $p_k$ to
be finite, and thus $\dim (H_k)=h-(k-4)$. Of course we are mainly interested
in smooth curves. We will say that (an irreducible component of) $H_k$ or $%
H_k^s$ has the {\it expected dimension} if it is of dimension $h-(k-4). $ }
\tit
The following general statement will be quite useful:
\tit
{\bf Proposition 2.4.}
{\it With notations as above, let $H_k',H_{k+1}'$ 
be irreducible components of $H_k,H_{k+1}$ such that 
$H_{k+1}^{^{\prime}}\subset H_k'.$ Set $\i _k'=p_k^{-1}(H_k').$ 
Assume $\dim (\i _k^{^{\prime}})=h-(k-4)$ 
and 
$\dim (H_{k+1}')=h-(k-3)$ 
(i.e. $\i _k'$ 
and $H_{k+1}'$ are both of the expected
dimension). Furthermore assume $\i _k'$ smooth. Then $H_k'$ 
also has the expected dimension: $\dim (H_k')=h-(k-4).$ }
\tit
{\bf Proof.} Since $H_{k+1}'\subset H_k',$ $\dim
(H_k')\geq \dim (H_{k+1}')=h-(k-3).$ Since 
$\dim (\i _k')=h-(k-4),$ if 
$H_k'=p_k(\i _k')$ 
has not the expected dimension then 
$\dim (H_k')=\dim (H_{k+1}'),$ 
and since they are both irreducible, $H_k'=H_{k+1}'.$
Moreover, for general $C$ in 
$H_k',p_k^{-1}(C)$ has dimension one. By generic smoothness we
may assume that the general fiber $p_k^{-1}(C)$ is a smooth, equidimensional
curve. On the other hand, $C\in $ $H_{k+1}',$ this means that
among the $k$-secants to $C$ there is at least one $(k+1)$-secant. This 
$(k+1)$-secant will count several time as a $k$-secant. in other words 
$p_k^{-1}(C)$ is singular at points $(Z,C)$ where $Z$ is on a $(k+1)$-secant
line to $C;$ this contradicts the smoothness of $p_k^{-1}(C)$ for general $C$
\tit
\vfill\eject
{\bf Complete intersections}
\tit
Abusing notations we will denote by $H(a,b)$ the Hilbert scheme of complete
intersections of type $(a,b),a\leq b.$ As it is well known $H(a,b)$ is
integral of dimension $h(a,b)$ where $h(a,b)=h^0({\cal O} _{{\bf P}}(a))
+h^0({\cal O} _{{\bf P}}(b))-h^0({\cal O} _{{\bf P}}(b-a))-2$, if $a<b,$ 
and where $h(a,b)=2.h^0({\cal O} _{{\bf P}}(a))-4$ if $a=b.$
\tit
{\bf Lemma 2.5.} {\it If $k\leq b$ the morphism $q_k:\i _k(a,b)\rightarrow Al^k$ is
surjective, smooth, of relative dimension $max\{h(a,b)-2k,h(a,b)-(a+1)-k\}.$ 
In particular if $k\leq a+1$ $\i _k(a,b)$ 
is smooth, irreducible, of dimension $h(a,b)-(k-4).$ }
\tit
{\bf Proof.} This follows essentially from the fact that every $Z$ in $Al^k$
is a complete intersection $(1,1,k)$ and hence
gives independant conditions to forms of degree $\geq k-1$.
\tit
{\bf Lemma 2.6.} {\it If $k=b$ and $a\geq 4$, $\i _b^s$ is non-empty and 
$p_b:\i _b^s\rightarrow H_b^s(a,b)$ 
is generically finite (in fact
birational if $a<b$). In particular $H^s_b(a,b)$ is irreducible
of dimension $h(a,b)-(a-3)$. }
\tit
{\bf Proof.} First assume $a<b.$
Take a line $L\subset {\bf P}^3.$ By [L] (see also Remark 1.3.c), if $F_a$ is a
sufficiently general surface of degree $a\geq 4,$ containing $L$ then $Pic(F_a)$ 
is generated by $L$ and the hyperplane section, in particular $F_a$
doesn't contain any further line. 
Indeed if $R$ is another line on $F_a$ then $R\sim cH+dL$. From $H\cdot R=1$ we get
$d=1-ca$. Since $L^2=2-a$, $R\cdot L=c+(ac-1)(a-2)$, but since
$0\leq R\cdot L\leq 1$, we get a contradiction.
We may assume $F_a$ smooth. Let $C=F_a\cap
F_b,$ $F_b$ a general surface of degree $b$. Then $C$ is smooth and $L$ is a
(proper) $b-$secant to $C.$ Now $C$ has no other $k$-secants, $k>a;$ indeed
such a secant would have to lie on $F_a.$

The case $a=b$ requires an extra argument.

Let $C$ be a smooth complete intersection of type $(b,b)$ and let $R$ be a $%
b $-secant to $C.$ Observe that $R$ is a proper $b$-secant to $C.$ Moreover
from the exact sequence 
$$
0\rightarrow \i _{C\cup R}\rightarrow \i _C \rightarrow 
{\cal O} _R (-b) \rightarrow 0 
$$
it follows that $h^0(\i _{C\cup R}(b))=1.$ So, to any $b$-secant, $%
R, $ there is associated one surface, $F_R,$ of the pencil ${\bf P}(H^0(%
\i _C(b))\simeq {\bf P}^1.$ If $C$ has infinite $b$-secants we get a
morphism from (a component of) the curve, $\Lambda ,$ of $b-$secants to $%
{\bf P}^1:\varphi :\Lambda \rightarrow {\bf P}^1:R\rightarrow F_R$. If $%
\varphi $ is constant then $C$ lies on a degree $b$ surface of $b-$secants
lines, otherwise $\varphi $ is surjective and {\it every} $F\in H^0(%
\i _C(b))$ contains a $b$-secant line to $C$. We will show that this
cannot happen if $C$ is sufficiently general.

Indeed take a general surface, $F_b,$ of degree $b\geq 4.$ We may assume
that $F_b$ doesn't contain any line (Noether-Lefschetz theorem). Now let $%
F_b^{^{\prime }}$ be a general surface of degree $b$ containing the line $L.$
We may assume $X=F_b\cap F_b^{^{\prime }}$ smooth. Clearly $L$ is a $b$%
-secant to $X.$ Assume that every fiber of $p_b:\i %
_b^s(b,b)\rightarrow H_b^s(b,b)$ has dimension one. By Lemma 2.5 and by
generic smoothness we may assume that the general fiber of $p_b$ is a
smooth, equidimensional curve. Since curves as $X$ above are clearly general
in $H_b^s(b,b),$ we may assume, with notations as above, that $L\in \Lambda $
(i.e. $L$ is not an isolated point of $p_b^{-1}(X)$). Then we get a
contradiction. Indeed, on one hand the morphism $\varphi $ has to be
constant because $F_b$ doesn't contain any line. On the other hand if $%
F=\lambda F_b+\mu F_b^{^{\prime }}$ is the surface of $b$-secants
corresponding to the point $\varphi (\Lambda ),$ then $L\subset F\cap
F_b^{^{\prime }}=X$, which is absurd (notice that $F\neq F_b^{^{\prime }}$
because by [L] we may assume that $L$ is the only line on $%
F_b^{^{\prime }}$)
\tit
{\bf Remark 2.7.} Observe that $H^s_k(a,b)=H^s_{a+1}(a,b)$ if 
$k\geq a+1$, indeed if $C= F_a\cap F_b$ then every $k$-secant to $C$
is contained in $F_a$ and thus is $b$-secant to $C$. 
\tit 
{\bf Corollary 2.8.} {\it 
For $a\geq 4$ and $4\leq k\leq b,H_k^s(a,b)$ is (non empty) integral, of
the expected dimension $h(a,b)-(k-4)$ if $4\leq k\leq a+1$ and of 
dimension $h(a,b)-(a-3)$ if $a+1\leq k \leq b$. }
\tit
{\bf Proof.} Since every $b$-secant is a $k$-secant for $k\leq b,$ by
Lemma 2.6 we get: $\i _k^s(a,b)\neq \emptyset .$ Clearly $%
\i _k^s(a,b)$ is open in $\i _k(a,b).$ By Lemma 2.5 we
conclude that $H_k^s(a,b)$ is irreducible. It remains to 
compute the dimension of
$H_k^s(a,b)$. For $a+1\leq k \leq b$ we combine Lemma 2.6 and Remark 2.7.
If $4\leq k\leq a+1$
we argue by descending induction on $k$ using Proposition 
2.4, the starting point being the case $k=a+1$ (combining Lemma 2.6 and
Remark 2.7). The
assumptions of Prop. 2.4 are satisfied since $H_{k+1}^s(a,b)\subset
H_k^s(a,b)$ and both are irreducible, and since $\i _k^s(a,b)$ is
smooth (Lemma 2.5).
\tit
{\bf Remark 2.9.}
If $a<4,$ our arguments break down, but a detailed analysis of curves
on low degree surfaces should give a precise description of $H_k^s(a,b)$. 
\tit
{\bf Rational curves}
\tit
The starting point is the following basic remark:
\tit
{\bf Remark 2.10.}
Let $\delta $ be a $\infty ^3$ linear system of degree $d$
divisors on ${\bf P}^1.$ Assume $\delta $ very ample so that $\delta $
yields an immersion $f:{\bf P}^1\rightarrow {\bf P}^3.$ The points $%
f(x_1),...,f(x_k)$ of ${\bf P}^3$ will be aligned if and only if $\delta $
contains a pencil having $x_1,...,x_k$ in its base locus.
\tit
Effective divisors of degree $d$ on ${\bf P}^1$ are parametrized by 
${\bf P}:={\bf P}({\cal O} _{{\bf P}^1} (d))\simeq {\bf P}^d,$ so a pencil is a
line in ${\bf P},$ i.e. a point in $G(1,{\bf P}).$ We will denote by
${\cal P} _k\subset G(1,{\bf P})$, 
the locus of pencils with base locus of
length $\geq k.$
\tit
{\bf Lemma 2.11.} {\it With notations as above, ${\cal P}_k$ is irreducible, of
dimension $2d-k-2.$ }
\tit
{\bf Proof.} Fix an effective divisor of degree $k$ on
${\bf P}^1:D=x_1+...+x_k.$ 
The effective divisors of degree $d$ containing $D$ build
a ${\bf P}^{d-k}$ in ${\bf P},$ let's denote it by ${\bf P}^{d-k}(D).$ A
line in ${\bf P}^{d-k}(D)$ is a pencil having $x_1,...,x_k$ in its base
locus. Now ${\cal P}_k$ can be described as follows: it is the image in 
$G(1,{\bf P})$ 
of a fibration over $Hilb^k{\bf P}^1\simeq {\bf P}^k$ with
fibers isomorphic to $G(1,d-k).$ Let's try to be more precise. Consider the
natural exact sequence:

$$
0\rightarrow K \rightarrow H^0({\cal O} _{{\bf P}^1}(k))\otimes 
H^0({\cal O} _{{\bf P}^1}(d-k))\rightarrow H^0({\cal O} _{{\bf P}^1}(d)) \rightarrow 0 
$$

We have the Segre embedding: ${\bf P}^k\times {\bf P}^{d-k}\hookrightarrow
S\subset {\bf P}^N,N=k(d-k)+d,$ where 
${\bf P}^N\simeq {\bf P(}H^0({\cal O}_{{\bf P}^1}(k))\otimes H^0({\cal O}_{{\bf P}^1}(d-k))),$ 
the projection
of $S$ from ${\bf P}(K)$ yields a finite morphism $S\rightarrow {\bf P}$
which presents ${\bf P}$ as ruled by the ${\bf P}^{d-k}$-fibers of $S.$ This
defines a diagram:

$$
\matrix{
S &  &  \cr 
\downarrow &  &  \cr 
{\bf P}^k & \rightarrow & Z\subset G(d-k,{\bf P}) \cr }
$$

where $Z$ is irreducible of dimension $k.$ Now in 
$G(1,{\bf P})\times G(d-k,{\bf P})$ 
consider the incidence variety: 
$P:=\left\{ (L,E)\in G(1,{\bf P)\times }Z/L\subset E\right\} ,$ 
clearly $P$ is irreducible, of dimension 
$k+\dim G(1,d-k).$ 
Finally ${\cal P }_k=\phi (P),$ where 
$\phi :P\rightarrow G(1,{\bf P})$ 
is the projection. To conclude we observe that $\phi $ 
is birational: let $L\subset {\bf P}^{d-k}(D),L$ generated by the
two points $D+D_{d-k},D+D_{d-k}^{\prime },$ where $D$ is a given effective
divisor of degree $k;$ if $Supp(D)\cap Supp(D_{d-k})=Supp(D)\cap
Supp(D_{d-k}^{\prime })=$ $Supp(D_{d-k})\cap Supp(D_{d-k}^{\prime
})=\emptyset ,$ then $L$ cannot be contained in ${\bf P}^{d-k}(D^{\prime })$
for $D\not = D^{\prime }$

\tit
We recall the following classical fact:
\tit
{\bf Lemma 2.12. }{\it For $d\geq 3$ there exists a smooth rational curve 
$C\subset {\bf P}^3$ 
with a $(d-1)-$secant line. Moreover if $C$ has more than one $(d-1)-$
secant lines, then $C$ lies on a smooth quadric surface. Hence for $d\geq 5,$
there exists a smooth rational curve of degree $d$ in ${\bf P}^3$ with
exactly one $(d-1)-$secant line.}
\tit
{\bf Proof.} We may assume $d\geq 5.$ For the existence look at rational
curves on the cubic scroll $S\subset {\bf P}^4.$ Then project from a general
point. The image $\overline{S}$ of $S$ is a cubic surface with a double
line, the double line will be the $(d-1)-$secant line. Finally observe that
if $C\subset \overline{S}$ then $h^0({\cal I} _C (2))=0$ . Indeed, if $d>6$
this is just by degree reasons, if $5\leq d\leq 6$ and if 
$h^0({\cal I} _C (2))\not = 0,$ 
then $C$ would have to be projectively normal (because it
would have to be a complete intersection if $d=6,$ or linked to a line if 
$d=5$),  
and this is absurd.

Finally if $R,L$ are two $(d-1)-$secant lines to $C,$ then 
$R\cap L=\emptyset $ 
(otherwise the plane $<R,L>$ will intersect $C$ in too many
points). Let $D$ be a three secant line to $C,$ $D$ is disjoint from $R$ and 
$L.$ Let $Q$ be the quadric generated by $R,L$ and $D$, then $Q$ intersects $C$ 
in $2d+1$ points, hence $C\subset Q$.
\tit
{\bf Lemma 2.13.}
{\it For $4\leq k\leq d-1,$ $H_k^s(d,0),$ the locus of smooth rational
curves of degree $d$ with a $k-$secant line, is irreducible. Moreover $%
H_{d-1}^s(d,0)$ has the expected dimension $3d+5.$}
\tit
{\bf Proof. }First of all notice that, by Lemma 2.12, 
$H_k^s(d,0)\not = \emptyset .$ 
A smooth rational curve of degree $d$ in ${\bf P}^3$ 
with a $k-$ secant line corresponds to a $\infty ^3-$linear system
containing a pencil having $k$ points in its base locus (see Remark 2.10). 
In $G(1,{\bf P})\times G(3,{\bf P)}$, consider the
(restricted) incidence variety: $I\subset {\cal P} _k\times G( 3,{\bf P}),$ 
$I=\left\{ (L,E)/L\subset E\right\} $ and the associated diagram: 
${\cal P} _k\buildrel q\over \longleftarrow I \buildrel p \over \longrightarrow G(3,{\bf P}).$
Since the fibers of $q$ are isomorphic to $G(1,d-2),$ by Lemma 2.11, we
conclude that $I$ is irreducible, of dimension $4d-k-8.$ It follows that $p(I)$ 
is irreducible. As noticed at the beginning of the proof, the
intersection, in $G(3,{\bf P}),$ of $p(I)$ with the open set of very ample $\infty ^3$ 
linear systems is non empty. Let $U$ denote this intersection.
Choosing a basis in the $4-$vector space corresponding to an $\infty ^3-$%
linear system, yields a Stiefel fibration ${\cal S}\rightarrow U$, with
fibers isomorphic to $Aut({\bf P}^3),$ now $Aut({\bf P}^1)$ acts (fiberwise)
on ${\cal S}$ and the quotient is $H_k^s(d,0).$ It follows that $H_k^s(d,0)$ 
is irreducible. From this description it also follows that 
$\dim (H_k^s(d,0))=\dim (p(I))+12,$ 
in particular $H_k^s(d,0)$ has the expected
dimension if and only if $p$ is generically finite, which amounts to say
that the generic rational curve of degree $d$ with a $k-$secant line has
only a finite number of $k-$secant lines. If $k=d-1,$ from Lemma 2.12, $p$
is birational and we are done (for $k<d-1$ we didn't find a simple direct
argument).
\tit
In the next lemma notations are as in 2.1.
\tit
{\bf Lemma 2.14.} {\it For $d\geq 5$ and $4\leq k\leq d-1,$ ${\cal I}_k^s$ is
smooth, of the expected dimension $4d-(k-4).$}
\tit
{\bf Proof.} Consider the incidence variety 
$${\bf I} _k=\left\{
(Z,C)\in Hilb^k{\bf P}^3\times H^s(d,0)/Z\subset C\right\} $$
and the
corresponding diagram:

$Hilb^k{\bf P}^3\buildrel f \over \longleftarrow  {\bf I} _k 
\buildrel \pi \over \longrightarrow H^s(d,0).$ 
Since $H^s(d,0)$ is smooth at $\left[ C\right] ,$
if the natural map $r:H^0(N_C)\rightarrow H^0(Z,N_C\mid Z)$ is surjective
then $f$ is smooth at $(Z,C)$ (see [K] and also [Pe]). Since 
$N_C\simeq {\cal O}_{{\bf P}^1}(a)\oplus   {\cal O}_{{\bf P}^1}(b),$ 
with
$d+2\leq a$, $d+2\leq b$, we have $h^1(N_C(-Z))=0$ and $r$ is surjective. We conclude that $f$
is smooth. Since ${\bf I}_k$ is irreducible of dimension $4d+k$, and
since $f$ is dominant [Pe], we conclude that $f$ has relative
dimension $4d-2k.$ Now ${\cal I} _k^s$ is obtained by base-change:

$$
\matrix{
{\cal I} _k^s & \rightarrow & {\bf I}_k \cr 
q\downarrow &  & \downarrow f \cr 
Al^k & \hookrightarrow & Hilb^k{\bf P}^3 \cr }
$$

Since ${\cal I} _k^s\not = \emptyset $ (Lemma 2.12), it follows that $q$
is smooth, of relative dimension $4d-2k.$ Observe that $q$ is surjective
(take a smooth quadric, $Q,$ containing the line $<Z>$ and look at the
linear system of rational curves of degree $d$ on $Q$ passing through $Z).$
So, since $Al^k$ is smooth, of dimension $4+k,$ ${\cal I }_k^s$ is smooth,
of dimension $4d-(k-4)$.
\tit 
{\bf Theorem 2.15.}
{\it For $d\geq 5$ and $4\leq k\leq d-1,$ $H_k^s(d,0)$ is integral, of the
expected dimension $4d-(k-4).$ Moreover if $C$ corresponds to a general
point of $H_k^s(d,0)$ then $C$ has a finite number, $s(k),$ of $k-$secant
lines (which are all proper and maximal).}
\tit
{\bf Proof. }
The proof is by descending induction on $k.$ The case $k=d-1$
follows from Lemma 2.12 and Lemma 2.13. Assume the theorem for $k+1 $ 
and suppose $\dim (H_k^s(d,0))<4d-(k-4).$ Since
$H_{k+1}^s(d,0)\subset H_k^s(d,0)$ 
and since they are both irreducible, it follows that
$H_k^s(d,0)=H_{k+1}^s(d,0),$ 
furthermore, for any irreducible component, 
$\overline{{\cal I}}_k^s$ 
of ${\cal I}_k^s,$ 
$\overline{{\cal I}} _k^s\rightarrow H_k^s(d,0)$ 
has one-dimensional fibers. Since 
${\cal I} _k^s$ 
is smooth (Lemma 2.14), we can repeat the argument of Prop.
2.4 and get a contradiction. This proves that $H_k^s(d,0)$ has the
expected dimension. Hence $H_{k+1}^s(d,0)$ is a proper closed subset of 
$H_k^s(d,0).$ 
So if $C$ corresponds to a general point of $H_k^s(d,0),C$ has
no $(k+1)-$secant lines. Moreover for every irreducible component 
$\overline{{\cal I}}_k^s$ of ${\cal I}_k^s$, 
$\overline{{\cal I}} _k^s\rightarrow H_k^s(d,0)$ 
is generically finite. Since there are a finite
number of irreducible components, $C$ has only finitely many $k-$secant lines.
\tit
{\bf Remarks 2.16.}
(i) We recover the classicaly known fact that every rational quintic curve
has a $4-$secant line.

(ii) We have $s(d-1)=1$ (note the analogy with the 
complete intersections
case), while $s(4)$ may be recover from the formula 
for quadrisecant lines.
As far as we know it is still an open problem to 
determine $s(k)$ in general.

\tit
\centerline{\bf References}
\tit
\item{[B]} Basili, B.: {\it Indice de Clifford des intersections 
compl$\grave e$tes.} Bull. Soc. Math. France {\bf 124}, 61-95 (1996).

\item{[CM]} Coppens, M.- Martens, G.: {\it Secant spaces and 
Clifford's Theorem.} Compositio Math. {\bf 78}, 193-212 (1991).

\item{[ELMS]} Eisenbud, D.-Lange, H.-Martens, G.-Schreyer F.O.:
{\it The Clifford dimension of a projective curve. } 
Compositio Math. {\bf 72}, 173-204 (1989).

\item{[H]} Hartshorne, R.: {\it Algebraic Geometry.} G.T.M., {\bf 52} (1977).

\item{[K]} Kleppe, J.O.: {\it Non reduced components of the Hilbert scheme of smooth space curves.}  "Space curves", Proc. Rocca di Papa 1985, L.N.M. {\bf 1266}, 181-207 (1987).

\item{[La]} Lazarsfeld, R.: {\it Lectures on linear series.} Alg-Geom /9408011. 

\item{[L]}  Lopez, A.F.: {\it Noether-Lefschetz theory and
Picard group of projective surfaces. } Memoirs of the AMS {\bf 438}
(1991).

\item{[P]} Paoletti, R.: {\it On Halphen's speciality theorem.} "Higher dimensional complex varieties", Proc. Trento, 1994 (de Gruyter ed.), 341-355 (1996). 

\item{[Pe]} Perrin, D.: {\it Courbes passant par m points g\'en\'eraux de $\Pt $.} Bull. Soc. Math. France, M\'emoire (nouv. s\'erie), {\bf 28-29} (1987).

\item{[T]}  Tyurin, A.N.: {\it Cycles, curves and vector bundles on an algebraic surface.} Duke Math. J., {\bf 54} , 1-26 (1987).

\end